\documentclass[12pt,a4paper]{article}

\newtheorem{theorem}{Theorem}[section]
\newtheorem{proposition}[theorem]{Proposition}
\newtheorem{lemma}[theorem]{Lemma}
\newtheorem{corollary}[theorem]{Corollary}

\usepackage{hyperref}
\usepackage{amsmath}
\usepackage{amssymb}
\usepackage{array}
\usepackage[english]{babel}
\usepackage{calc}
\usepackage{enumerate}
\usepackage[latin1]{inputenc}
\usepackage{latexsym}

\begin{document}
\title{A class of weighted convolution Fr\'echet algebras}
\author{Thomas Vils Pedersen}

\maketitle


\begin{abstract}
\noindent
For an increasing sequence $(\omega_n)$ of algebra weights on $\mathbb R^+$ we
study various properties of the Fr\'echet algebra 
$A(\omega)=\bigcap_n L^1(\omega_n)$ obtained
as the intersection of the weighted Banach algebras $L^1(\omega_n)$. 
We show that every endomorphism of $A(\omega)$ is standard, if for all
$n\in\mathbb N$ there exists $m\in\mathbb N$ such that
$\omega_m(t)/\omega_n(t)\to\infty$ as $t\to\infty$.
Moreover, we characterise the continuous derivations on this algebra: 
If for all $n\in\mathbb N$ there exists $m\in\mathbb N$ such that 
$t*\omega_n(t)/\omega_m(t)$ is bounded on $\mathbb R^+$, then the continuous  
derivations on $A(\omega)$ are exactly the linear maps $D$ of the form 
$D(f)=(Xf)*\mu$ for $f\in A(\omega)$, where $\mu$ is a measure in 
$B(\omega)=\bigcap_nM(\omega_n)$ and $(Xf)(t)=tf(t)$ for $t\in\mathbb R^+$ and
$f\in A(\omega)$. If the condition is not satisfied, we show that
$A(\omega)$ has no non-zero derivations. 

\end{abstract}

\section{Introduction}
\label{sec:intro}

In this paper we will study a class of Fr\'{e}chet algebras, which are
intersections
(or more formally projective limits) of decreasing sequences of
weighted convolution Banach algebras on the half-line. 
We start by recalling the definition and some basic properties of
these algebras (see for instance \cite{Gr} for further details).

Let $L^1(\mathbb R^+)$ be the Banach space of (equivalence classes of)
integrable functions $f$ on $\mathbb R^+=[0,\infty)$ with the norm 
$\|f\|=\int_0^{\infty}|f(t)|\,dt$.
A positive Borel function $\omega$ on $\mathbb R^+$ is called an 
{\it algebra weight} if 
(a) $\omega$ and $1/\omega$ are locally bounded on $\mathbb R^+$,
(b) $\omega$ is right continuous on $\mathbb R^+$,
(c) $\omega$ is submultiplicative, that is $\omega(t+s)\le\omega(t)\omega(s)$
for $t,s\in\mathbb R^+$, and (d) $\omega(0)=1$.
We then define $L^1({\omega})$ as the weighted space of functions $f$ on 
$\mathbb R^+$ 
for which $f\omega\in L^1(\mathbb R^+)$ with the inherited norm 
$$\|f\|_{\omega}=\int_0^{\infty}|f(t)|\omega(t)\,dt.$$
With the usual convolution product
$$(f\ast g)(t)=\int_0^t f(s)g(t-s)\,ds\qquad
\text{for $t\in\mathbb R^+$ and $f,g\in L^1({\omega})$}$$
it is well known that $L^1({\omega})$ is a commutative Banach algebra.
Similarly, the space $M(\omega)$ of locally finite complex Borel measures 
$\mu$ on $\mathbb R^+$ for which 
$$\|\mu\|_{\omega}=\int_0^{\infty}\omega(t)\,d|\mu|t<\infty$$
is a Banach algebra under convolution and contains $L^1({\omega})$ as a closed 
ideal.

We are now ready to define the Fr\'{e}chet algebras which we will
study.  Throughout
this paper $\omega=(\omega_n)$ will denote an increasing sequence of
algebra weights on $\mathbb R^+$ and we let
$$A(\omega)=\bigcap_n L^1(\omega_n)\qquad\text{and}\qquad 
B(\omega)=\bigcap_nM(\omega_n).$$
Equipped with the increasing sequence of norms 
$\|\mu\|_n=\|\mu\|_{\omega_n}\ (\mu\in B(\omega))$, it is easily seen
that $A(\omega)$ and $B(\omega)$ become Fr\'{e}chet algebras. 
More abstractly, these Fr\'{e}chet algebras can be viewed as
projective limits of
the weighted Banach algebras in question (with the inclusion maps).
In particular, sets of the form 
$\{g\in A(\omega):\|g\|_{L^1(\omega_m)}<\delta\}$ with
$\delta>0$ and $m\in\mathbb N$ form a base at zero for the projective limit
topology on $A(\omega)$. For a sequence $(f_k)$ in $A(\omega)$ we thus have
$f_k\to0$ in $A(\omega)$ as $k\to\infty$ if and only if $f_k\to0$ in
$L^1(\omega_n)$ as $k\to\infty$ for every $n\in\mathbb N$.
We refer to Michael's original memoir \cite{Mic} for general
background material on Fr\'{e}chet algebras.

In the rest of the paper we will make the following assumptions on the
weights $\omega=(\omega_n)$:
\begin{enumerate}[(a)]
\item 
$\omega_n(t)\to\infty$ as $t\to\infty$ for every $n\in\mathbb N$,
\item
$\lim_{t\to\infty}\omega_n(t)^{1/t}=1$ for every $n\in\mathbb N$,
\item
$\sup_{t\in\mathbb R^+}\omega_{n+1}(t)/\omega_n(t)=\infty$ for every  
$n\in\mathbb N$.
\end{enumerate}
It follows from (a) that the weights are semisimple, and the
assumption (b) is equivalent to each of the algebras $L^1(\omega_n)$ having
the right hand half-plane as character space
(see the proof of Theorem~\ref{th:char}). 
Moreover, condition (c) means that $L^1(\omega_{n+1})\subset L^1(\omega_n)$, so
that the descending chain $(L^1(\omega_n))$ does not stabilise.
We are mainly interested in the case where we further have
\begin{enumerate}
\item[(d)] 
$\omega_n(t)\to\infty$ as $n\to\infty$ for every $t\in\mathbb R^+$, 
\end{enumerate}
but we will
see that there are also interesting aspects in the case where the 
``limit weight'' 
$\omega_{\infty}(t)=\lim_{n\to\infty}\omega_n(t)$ is finite 
(for instance $\omega_n(t)=(1+t)^{1-1/n}$ for $t\in\mathbb R^+$ and 
$n\in\mathbb N$).
It may also be of interest to study $A(\omega)$ for radical weights
$(\omega_n)$ with $\omega_{\infty}$ either radical or semisimple.

The organisation of the paper is as follows: In
Section~\ref{sec:various} we collect various basic results about
the algebra $A(\omega)$ including a characterisation of $B(\omega)$ as
the multiplier
algebra of $A(\omega)$. Endomorphisms of $A(\omega)$ are studied in 
Section~\ref{sec:endo}. In particular, under a slightly stronger
assumption than (c) above we show that every endomorphisms of $A(\omega)$
is standard. Finally, in Section~\ref{sec:der} we describe the 
derivations on $A(\omega)$.

\section{Basic results}
\label{sec:various}

In this section we present som basic results about the algebras $A(\omega)$
and $B(\omega)$. The following result is a direct consequence of the
topology on $A(\omega)$.

\begin{lemma}
\label{le:top}
A linear map $T:A(\omega)\to A(\omega)$ is continuous if and only if
for all $n\in\mathbb N$
there exists $m\in\mathbb N$ such that $T$ extends continuously to a
map (also denoted) 
$T:L^1(\omega_m)\to L^1(\omega_n)$.

Similarly, a linear functional $\varphi:A(\omega)\to\mathbb C$ is
continuous if and only if 
$\varphi$ extends continuously to $L^1(\omega_n)$ for some $n\in\mathbb N$.
\end{lemma}

\noindent{\bf Proof}\quad
Clearly $T:A(\omega)\to A(\omega)$ is continuous if the continuous extensions
exist. Conversely, assume that $T:A(\omega)\to A(\omega)$ is continuous and let
$n\in\mathbb N$. Since $U=\{f\in A(\omega):\|f\|_n<1\}$ is an open
neighbourhood of 0 in $A(\omega)$, the same holds for $T^{-1}(U)$,
so there exist $\delta>0$ and $m\in\mathbb N$ such that 
$\{g\in A(\omega):\|g\|_{L^1(\omega_m)}<\delta\}\subseteq T^{-1}(U)$. Hence
$\|Tg\|_n\le\frac{1}{\delta}\,\|g\|_m$ for $g\in A(\omega)$, so a standard
argument using Cauchy sequences shows that $T$ extends continuously to a
map $T:L^1(\omega_m)\to L^1(\omega_n)$. Similarly for the second part
of the lemma.
{\nopagebreak\hfill\raggedleft$\Box$\bigskip}

For a weight $\omega$ let $L^{\infty}(1/\omega)$ denote the the Banach
space of functions $h$ on $\mathbb R^+$ for which 
$\|h\|=\sup_{t\in\mathbb R^+}|h(t)|/\omega(t)<\infty$. It is well known that
the duality
$\langle f,h\rangle=\int_0^{\infty}f(t)h(t)\,dt$ for 
$f\in L^1(\omega)$ and $h\in L^{\infty}(1/\omega)$ identifies 
$L^{\infty}(1/\omega)$
isometrically isomorphic with the dual space of $L^1(\omega)$. 
As a consequence of the previous lemma
we therefore have the following characterisation of the dual space of
$A(\omega)$. 

\begin{corollary}
\label{cor:dual}
Every function $h\in\bigcup_n L^{\infty}(1/\omega_n)$ defines a continuous
linear functional on $A(\omega)$ by 
$\langle f,h\rangle=\int_0^{\infty}f(t)h(t)\,dt\ (f\in A(\omega))$
and conversely every continuous linear functional on $A(\omega)$ is of
this form.
\end{corollary}

\begin{theorem}
\label{th:char}
Every character on $A(\omega)$ is continuous. Moreover, for every  
$z\in\mathbb C$
with $\text{Re}\,z\ge0$ the Laplace transform
$${\cal L}(f)(z)=\int_0^{\infty}f(t)e^{-zt}\,dt\qquad(f\in A(\omega))$$
defines a continuous character on $A(\omega)$ and conversely every
character on $A(\omega)$
has this form.
\end{theorem}

\noindent{\bf Proof}\quad
For $z\in\mathbb C$ with $\text{Re}\,z\ge0$ the map 
$f\mapsto{\cal L}(f)(z)\ (f\in A(\omega))$ defines a continuous
character on each
$L^1(\omega_n)$ by the usual characterisation of the characters on 
$L^1(\omega_n)$
(\cite[Theorem~4.7.27]{Da:Book}) and hence defines a continuous
character on $A(\omega)$. 
Let $g(t)=e^{-t}\ (t\ge0)$. It follows from
\cite[Theorem~4.7.26]{Da:Book} that the constant function with value 1
is a polynomial generator of $L^1(e^{-t}\omega_n(t))$ and thus that
$g$ is a polynomial generator of $L^1(\omega_n)$ for every
$n\in\mathbb N$. 
Hence $g$
is a polynomial generator of $A(\omega)$, so every character on $A(\omega)$ is
continuous by \cite[Corollary~4.10.11]{Da:Book}. 
By Lemma~\ref{le:top} it thus follows that every character
$\varphi:A(\omega)\to\mathbb C$ extends continuously to $L^1(\omega_n)$ for some
$n\in\mathbb N$. Hence there exists $z\in\mathbb C$ with 
$\text{Re}\,z\ge0$ such that $\varphi(f)={\cal L}(f)(z)$ for
$f\in L^1(\omega_n)$.
{\nopagebreak\hfill\raggedleft$\Box$\bigskip}

It follows from Theorem~\ref{th:char} that the maximal modular ideals
in $A(\omega)$ are exactly the sets $\{f\in A(\omega):{\cal L}(f)(z)=0\}$,
where $z\in\mathbb C$ with $\text{Re}\,z\ge0$. The corresponding sets 
$\{\mu\in B(\omega):{\cal L}(\mu)(z)=0\}$ are maximal ideals in
$B(\omega)$ and their 
intersection equals $\{0\}$, so we have the following result.

\begin{corollary}
\label{co:ss}
The Fr\'{e}chet algebras $A(\omega)$ and $B(\omega)$ are semisimple.
\end{corollary}

It is well known that each $L^1(\omega_n)\ (n\in\mathbb N)$ has a
bounded approximate identity, for instance 
$e_k(t)=k\cdot1_{[0,1/k]}\ (k\in\mathbb N)$. 
Hence $e_k*f\to f$ in $A(\omega)$ as
$k\to\infty$ for $f\in A(\omega)$, so $(e_k)$ is also a bounded approximate
identity for $A(\omega)$. However, if $\lim_{n\to\infty}\omega_n(t)=\infty$
for every $t\in\mathbb R^+$
(which is the case we are most interested in), then $\|f\|_n\to\infty$
as $n\to\infty$ for every non-zero $f\in A(\omega)$. Hence $A(\omega)$ does not
possess a uniformly 
bounded approximate identity, that is, a bounded approximate identity
where the bound in $L^1(\omega_n)$ is independent of $n\in\mathbb N$. Note, however,
that $\|e_k\|_n\to1$ as $k\to\infty$ for every $n\in\mathbb N$.

For a Fr\'{e}chet algebra resp.\ Fr\'echet module with a uniformly bounded
approximate identity, Craw (\cite{Cr}) resp.\ Summers (\cite{Su})
generalised Cohen's factorisation theorem. (See also \cite{Do-Wi} for
a thorough discussion of approximate identities and factorization.) 
These results do not apply to $A(\omega)$ if  
$\lim_{n\to\infty}\omega_n(t)=\infty$ for every $t\in\mathbb R^+$ and we do
not know whether we have factorisation in $A(\omega)$ or in Fr\'echet  
$A(\omega)$-modules; not even in the simplest sense of being able to factor every
$f\in A(\omega)$ as a product $f=g*h$ with $g,h\in A(\omega)$. 
We complement the discussion by mentioning that in the special
case where there exists $a>0$ such that 
$\sup_{n\in\mathbb N,\,0\le t\le a}\omega_n(t)<\infty$, the
sequence $(e_k)$ {\it is} a uniformly bounded approximate identity for $A(\omega)$,
and we therefore have factorisation in $A(\omega)$.

A linear map $T$ on a commutative Fr\'{e}chet algebra $\cal B$ is called a multiplier if 
$T(ab)=T(a)b$ for every $a,b\in\cal B$. 
Recall from \cite[Theorem~2.2]{Gh:Hom} that every $\mu\in M(\omega)$ defines
a continuous multiplier $T_{\mu}$ on $L^1(\omega)$ by 
$T_{\mu}(f)=\mu\ast f\ (f\in L^1(\omega))$ and that this identifies $M(\omega)$  
isometrically isomorphic with the multiplier algebra of $L^1(\omega)$. For
the algebra $A(\omega)$ we have a similar result.

\begin{theorem}
\label{th:mul}
Every $\mu\in B(\omega)$ defines a continuous multiplier $T_{\mu}$ on $A(\omega)$ by 
$T_{\mu}(f)=\mu\ast f\ (f\in A(\omega))$ and every multiplier on $A(\omega)$ has
this form. In particular, every multiplier on $A(\omega)$ is continuous.
\end{theorem}

\noindent{\bf Proof}\quad
Clearly $T_{\mu}$ is a continuous multiplier on $A(\omega)$ for every
$\mu\in B(\omega)$. Conversely, let $T$ be a multiplier on $A(\omega)$. Since
$A(\omega)$ is an integral domain, it is well known 
(see \cite[Proposition~2.5.12]{Da:Book}) that $T$ is automatically
continuous.
Let $(e_k)$ be a bounded approximate identity for $A(\omega)$ and let
$n\in\mathbb N$. By Lemma~\ref{le:top}
there exists $m\in\mathbb N$ such that $T$ extends continuously to a map 
$T:L^1(\omega_m)\to L^1(\omega_n)$. Hence the sequence $(Te_k)$ is bounded in $L^1(\omega_n)$,
and it follows from the proof of \cite[Theorem~3.3.40]{Da:Book} that 
the sequence has a weak-star cluster point $\mu_n\in M(\omega_n)$ for which 
$T(f)=\mu_n\ast f$ for $f\in L^1(\omega_n)$. By uniqueness of $\mu_n$ it is
independent on $n$ and the result follows.
{\nopagebreak\hfill\raggedleft$\Box$\bigskip}

The identification in the previous theorem of $B(\omega)$ as the multiplier
algebra of $A(\omega)$ induces a strong operator
topology (denoted by SO) on $B(\omega)$ in which a net $(\mu_{\beta})$
tends to 0 if and only if $\mu_{\beta}\ast f\to0$ for every $f\in A(\omega)$.
We will now introduce another topology on $B(\omega)$, denoted by $\sigma$, 
which is the projective limit topology, when each of the spaces
$M(\omega_n)$ is regarded as the multiplier algebra of $L^1(\omega_n)$ and is equipped
with the corresponding strong operator topology (denoted by SO$_n$). 
Hence $\sigma$ is the weakest topology on $B(\omega)$ making all the inclusions 
$\iota_n:(B(\omega),\sigma)\to(M(\omega_n),\text{SO}_n)$ for $n\in\mathbb N$ continuous 
(see, for instance, \cite[Proof of Proposition~2.7]{Mic}).
Also, for a net $(\mu_{\beta})$ in $B(\omega)$ we have 
$\mu_{\beta}\to0$ in the $\sigma$-topology if and only if $\mu_{\beta}\to0$
strongly in $M(\omega_n)$ for every $n\in\mathbb N$. 
We will need the following results about the $\sigma$-topology.

\begin{lemma}
\label{le:SOsi}
We have SO\,$\subseteq\sigma$, that is, SO is weaker than $\sigma$.
\end{lemma}

\noindent{\bf Proof}\quad
Let $(\mu_{\beta})$ be a net in $B(\omega)$ with 
$\mu_{\beta}\to0$ in the $\sigma$-topology, that is, $\mu_{\beta}\to0$
strongly in $M(\omega_n)$ for every $n\in\mathbb N$.
For $f\in A(\omega)$ we thus have $\mu_{\beta}*f\to0$ in $L^1(\omega_n)$ for every
$n\in\mathbb N$ and thus $\mu_{\beta}*f\to0$ in $A(\omega)$. Hence 
$\mu_{\beta}\to0$ strongly in $B(\omega)$, so the identity map
$\iota:(B(\omega),\sigma)\to(B(\omega),\text{SO})$ is continuous and the result follows.
{\nopagebreak\hfill\raggedleft$\Box$\bigskip}

We do not know whether the $\sigma$- and the SO-topologies are
identical, but on bounded sets in $B(\omega)$ (that is, bounded in every $M(\omega_n)$) 
that is the case: Let $(\mu_{\beta})$ be a bounded net in $B(\omega)$ with 
$\mu_{\beta}\to0$ strongly and let $n\in\mathbb N$. Then 
$\mu_{\beta}*f\to0$ in $L^1(\omega_n)$ for every $f\in A(\omega)$ and by the
boundeness thus for $f\in L^1(\omega_n)$. Hence $\mu_{\beta}\to0$ strongly in
$M(\omega_n)$ for every $n\in\mathbb N$ and thus $\mu_{\beta}\to0$ in the $\sigma$-topology.
Moreover, we have the following result.

\begin{proposition}
\label{prop:dualSOsi}
The dual spaces of $(B(\omega),\sigma)$ and $(B(\omega),\text{SO})$ are the
same, that is, a linear functional on $B(\omega)$ is $\sigma$-continuous if and
only if it is strongly continuous.
\end{proposition}

\noindent{\bf Proof}\quad
It follows from Lemma~\ref{le:SOsi} that a strongly continuous linear functional
on $B(\omega)$ also is $\sigma$-continuous. Conversely, let $\varphi:B(\omega)\to\mathbb C$ 
be a $\sigma$-continuous linear functional. 
Since $\sigma$ is the weakest topology on $B(\omega)$ making all the inclusions 
$\iota_n:(B(\omega),\sigma)\to(M(\omega_n),\text{SO}_n)$ for $n\in\mathbb N$ continuous, 
it follows that the family of seminorms
$$p_{n,f}(\mu)=\|\mu*f\|_n\qquad(\mu\in B(\omega))$$
with $n\in\mathbb N$ and $f\in L^1(\omega_n)$ defines the $\sigma$-topology.  
Hence (see, for instance, \cite[Theorem~IV.3.1]{Co}) there exist
$n_1,\ldots,n_J$ and $f_j\in L^1({\omega_{n_j}})\ (j=1,\ldots,J)$ such that
$$\|\varphi(\mu)\|\le\sum_{j=1}^J\|\mu*f_j\|_{n_j}
\le\sum_{j=1}^J\|\mu*f_j\|_N\qquad(\mu\in B(\omega))$$
with $N=\max\{n_1,\ldots,n_J\}$.
A standard argument using Cauchy sequences now shows that $\varphi$
extends to a linear functional $\varphi_N:M({\omega_N})\to\mathbb C$ with
$\|\varphi_N(\mu)\|\le\sum_{j=1}^J\|\mu*f_j\|_N$ for $\mu\in
M({\omega_N})$.
In particular this shows that $\varphi$ is strongly continuous on $B(\omega)$.
{\nopagebreak\hfill\raggedleft$\Box$\bigskip}

The previous proposition does not in itself imply that the $\sigma$-
and the SO-topologies are identical. For instance, if $X$ is a Banach
space, then the dual spaces of $X$ with respect to the norm and the
weak topologies are identical (\cite[p.\,63]{Ru:Funct}), whereas the
norm and the weak topologies themselves are not.

\section{Endomorphisms}
\label{sec:endo}

Homomorphisms between weighted convolution Banach algebras have been
studied extensively; see for instance \cite{Gr}, \cite{Gh:Hom},
\cite{Gh-Mc-Gr} and \cite{Gh-Gr:St}. For endomorphisms of $A(\omega)$ we
start with the following consequence of Theorem~\ref{th:char}. 

\begin{corollary}
Every endomorphism of $A(\omega)$ is continuous. 
More generally, if $\cal B$ is a Fr\'{e}chet algebra on which every character is
continuous, then every homomorphism $\Phi:{\cal B}\to A(\omega)$ is continuous.
\end{corollary}

\noindent{\bf Proof}\quad
This follows from Theorem~\ref{th:char} as in the proof of  
\cite[Theorem~2.3.3]{Da:Book} since $A(\omega)$ is semisimple 
(Corollary~\ref{co:ss}) and since the closed graph theorem holds for
operators on Fr\'{e}chet spaces (\cite[Theorem~A.3.25]{Da:Book}).
{\nopagebreak\hfill\raggedleft$\Box$\bigskip}

The next result (which can be generalised to homomorphisms from $A(\omega)$
to $A(\widetilde{\omega})$) is similar to and follows rather easily
from  the corresponding result for $L^1(\omega)$
(\cite[Theorems~3.4 and 3.6]{Gr}).

\begin{theorem}
\label{th-hom}
Let $\Phi$ be a non-zero endomorphism of $A(\omega)$. 
Then $\Phi$ has a unique extension to a continuous endomorphism 
$\widetilde{\Phi}$ of $B(\omega)$. 
Also, $\nu^t=\widetilde{\Phi}(\delta_t)\ (t\in\mathbb R^+)$ defines a
semigroup in $B(\omega)$ which is strongly continuous for $t>0$. Moreover,
$$\widetilde{\Phi}(\mu)=\int_0^{\infty}\nu^t\,d\mu(t)
\qquad\text{for }\mu\in B(\omega),$$
where the integrals exist as strong Bochner integrals in $B(\omega)$, that is,  
$\widetilde{\Phi}(\mu)*f=\int_0^{\infty}\nu^t*f\,d\mu(t)$ exists as a
Bochner integral in each $L^1(\omega_n)\ (n\in\mathbb N)$ for $\mu\in B(\omega)$ and $f\in A(\omega)$.
\end{theorem}

\noindent{\bf Proof}\quad
By the previous theorem $\Phi$ is continuous. For $n\in\mathbb N$ it thus follows
from Lemma~\ref{le:top} that there exists $m\in\mathbb N$ such that $\Phi$
extends continuously and uniquely to a homomorphism
$\Phi_n:L^1(\omega_m)\to L^1(\omega_n)$. By \cite[Theorems~3.4 and 3.6]{Gr}, $\Phi_n$ extends 
uniquely to a continuous homomorphism
$\widetilde{\Phi}_n:M(\omega_m)\to M(\omega_n)$, the semigroup $(\nu^t)$ is
strongly continuous in $L^1(\omega_n)$ for $t>0$ and 
$\widetilde{\Phi}_n(\mu)*f=\int_0^{\infty}\nu^t*f\,d\mu(t)$ exists as a
Bochner integral in $L^1(\omega_n)$ for $\mu\in M(\omega_m)$ and $f\in L^1(\omega_n)$. 
By uniqueness $\widetilde{\Phi}(\mu)=\widetilde{\Phi}_n(\mu)$ is
independent of $n\in\mathbb N$ for $\mu\in B(\omega)$, and the result follows.
{\nopagebreak\hfill\raggedleft$\Box$\bigskip}

In the rest of the section we will investigate other continuity
properties of the endomorphism $\widetilde{\Phi}$ of $B(\omega)$; namely
with respect to a wk$^*$ topology and the strong topology on $B(\omega)$.
If we denote by $C_0(1/\omega)$ the closed subspace of
$L^{\infty}(1/\omega)$ consisting of continuous functions $h$ on
$\mathbb R^+$ for which $h/\omega$ is vanishes at infinity, then it is
well known that $M(\omega)$ is isometrically isomorphic to the dual space of 
$C_0(1/\omega)$ with the duality being defined by
$$\langle h,\mu\rangle
=\int_0^{\infty}h(t)\,d\mu(t)\qquad(h\in C_0(1/\omega),\mu\in M(\omega)).$$
We can use this to show that $B(\omega)$ is also a dual space. Let 
$$D(1/\omega)=\bigcup_{n\in\mathbb N}C_0(1/\omega_n)$$
and equip $D(1/\omega)$ with the inductive limit topology
(see, for instance, \cite[Chapter~IV.5]{Co}). Since a
linear functional on $D(1/\omega)$ is continuous if and only if its
restriction to each $C_0(1/\omega_n)$ is continuous 
(see, for instance, \cite[Proposition~IV.5.7]{Co}) we obtain the
following result.

\begin{proposition}
\label{pr-dual}
The duality 
$$\langle h,\mu\rangle
=\int_0^{\infty}h(t)\,d\mu(t)\qquad(h\in D(1/\omega),\mu\in B(\omega))$$
identifies $B(\omega)$ with the dual space of $D(1/\omega)$. Moreover, for a
net $(\mu_{\beta})$ in $B(\omega)$ we have $\mu_{\beta}\to0$ wk$^*$ in
$B(\omega)$ if and only if $\mu_{\beta}\to0$ wk$^*$ in $M(\omega_n)$ for every
$n\in\mathbb N$.
\end{proposition}

For a homomorphism $\Phi:L^1(\omega_1)\to L^1(\omega_2)$ Grabiner
(\cite[Theorem~1.1]{Gr-wks}) proved that the extension
$\widetilde{\Phi}:M(\omega_1)\to M(\omega_2)$ is automatically wk$^*$
continuous. (See \cite{Pe:wkprop} for related
results about homomorphisms from $L^1(\omega)$ into other Banach algebras.) 
The corresponding result for endomorphisms of $A(\omega)$ is an almost
direct consequence of Grabiner's result.

\begin{theorem}
\label{th-wkscts}
Let $\Phi$ be an endomorphism of $A(\omega)$, let $\widetilde{\Phi}$ be the
unique extension to an endomorphism of $B(\omega)$ and let $(\nu^t)$ be 
the semigroup in $B(\omega)$ given by Theorem~\ref{th-hom}. 
Then $\widetilde{\Phi}$ is wk$^*$ continuous. Moreover, $(\nu^t)$ is
wk$^*$ continuous in $B(\omega)$ for $t\ge0$ and for every $n\in\mathbb N$ there exists
$m\in\mathbb N$ such that $\nu^t\omega_n(t)/\omega_m(t)\to0$ wk$^*$
in $M(\mathbb R^+)$ as $t\to\infty$.   
\end{theorem}

\noindent{\bf Proof}\quad
Let $(\mu_{\beta})$ be a net in $B(\omega)$ with $\mu_{\beta}\to0$ wk$^*$ in
$B(\omega)$ and let $n\in\mathbb N$. By Proposition~\ref{pr-dual} we have 
$\mu_{\beta}\to0$ wk$^*$ in $M(\omega_n)$. It follows from the proof of 
Theorem~\ref{th-hom} that there exists $m\in\mathbb N$ such that 
$\widetilde{\Phi}$ extends to a continuous homomorphism 
$\widetilde{\Phi}_n:M(\omega_m)\to M(\omega_n)$. By Grabiner's result
$\widetilde{\Phi}_n$ is wk$^*$ continuous, so 
$\widetilde{\Phi}(\mu_{\beta})=\widetilde{\Phi}_n(\mu_{\beta})\to0$
wk$^*$ in $M(\omega_n)$. Hence $\widetilde{\Phi}(\mu_{\beta})\to0$ wk$^*$ in
$B(\omega)$, so we deduce that $\widetilde{\Phi}$ is wk$^*$ continuous. 
Since $\delta_t$ is wk$^*$ continuous in $B(\omega)$ for $t\ge0$, it follows that
$\nu^t=\widetilde{\Phi}(\delta_t)$ is wk$^*$ continuous in $B(\omega)$ for $t\ge0$.
Similarly $\delta_t/\omega_m(t)\to0$ wk$^*$ in $M(\omega_m)$ as
$t\to\infty$, so 
$\nu^t/\omega_m(t)=\widetilde{\Phi}(\delta_t/\omega_m(t))\to0$ wk$^*$
in $M(\omega_n)$, that is, $\nu^t\omega_n(t)/\omega_m(t)\to0$ wk$^*$
in $M(\mathbb R^+)$ as $t\to\infty$. 
{\nopagebreak\hfill\raggedleft$\Box$\bigskip}

For homomorphisms $\Phi:L^1(\omega_1)\to L^1(\omega_2)$ it is not known
whether the semigroup $(\nu^t)$ is strongly continuous in $M(\omega_2)$ at
$t=0$. In the papers \cite{Gh-Gr:St} and \cite{Gh-Gr:Conv} this
problem was linked to the notion of {\it convergence factors}. Using
the wk$^*$ continuity of $(\nu^t)$ at $t=0$ and a result from \cite{Gh-Gr:Conv}
we can show that under a rather mild growth condition on the weights  
$\omega_n$ as $n\to\infty$, the semigroup $(\nu^t)$ from  
Theorem~\ref{th-hom} is strongly continuous in $B(\omega)$ at $t=0$.

\begin{theorem}
\label{th-sgcts}
Let $\Phi$ be a non-zero endomorphism of $A(\omega)$, and let $(\nu^t)$ be 
the semigroup from Theorem~\ref{th-hom}. Suppose that 
\begin{equation}
\label{eq:wein}
\text{for every $n\in\mathbb N$ there exists $m\in\mathbb N$ such that}\quad
\frac{\omega_m(s)}{\omega_n(s)}\to\infty\quad\text{as }s\to\infty.
\end{equation}
Then $(\nu^t)$ is strongly continuous in $B(\omega)$ for $t\ge0$.
\end{theorem}

\noindent{\bf Proof}\quad
By  Theorem~\ref{th-hom} we only need to prove the strong continuity
at $t=0$. Let $n\in\mathbb N$ and choose $m\in\mathbb N$ such that
$\omega_n(s)/\omega_m(s)\to0$ as $s\to\infty$. Let
$\eta=\omega_n/\omega_m$. Then $\eta$ is bounded and 
$M(\omega_m\eta)=M(\omega_n)$ is translation invariant. Moreover, 
$$\frac{\omega_m(r+s)}{\omega_m(s)}\,\eta(r+s)
=\frac{\omega_n(r+s)}{\omega_m(s)}
\le\omega_n(r)\,\frac{\omega_n(s)}{\omega_m(s)}\to0
\qquad\text{as }s\to\infty$$
for every $r>0$. It thus follows from \cite[Theorem~3.2]{Gh-Gr:Conv}
that $\eta$ is a convergence factor for $\omega_m$ at $0$. This means
that if $(\mu_k)$ is a sequence in $M(\omega_m)$ with $\mu_k\to0$
wk$^*$ in $M(\omega_m)$ as $k\to\infty$, then $\mu_k*f\to0$ in the norm of 
$L^1(\omega_m\eta)=L^1(\omega_n)$ as $k\to\infty$ for every
$f\in L^1(\omega_n)$, that is, $\mu_k\to0$ strongly in $M(\omega_n)$ as $k\to\infty$. 
By Theorem~\ref{th-wkscts} we have  
$\nu^t\to\widetilde{\Phi}(\delta_0)=\delta_0$ wk$^*$ in 
$M(\omega_m)$ as $t\to0$, so we conclude that $\nu^t\to\delta_0=\nu^0$
strongly in $M(\omega_n)$ and thus in $B(\omega)$ as $t\to0$.
{\nopagebreak\hfill\raggedleft$\Box$\bigskip}

We observe that condition (\ref{eq:wein}) is only a slight strengthening
of the standing assumption 
(c): $\sup_{t\in\mathbb R^+}\omega_{n+1}(t)/\omega_n(t)=\infty$ for every
$n\in\mathbb N$; see the introduction. In order to show that (\ref{eq:wein})
is a strictly stronger condition than (c), we construct below an
unbounded weight
$\omega$ which does not satisfy $\lim_{t\to\infty}\omega(t)=\infty$.
We then let $\omega_n=\omega^n$ and observe that $(\omega_n)$
satisfies (c) but not (\ref{eq:wein}). 
The weight $\omega(t)=2^{v(t)}$ is constructed as follows. 
For $m\in\mathbb N$ let $v(m)$ be the number of 1's in the binary expansion
of $m$ or equivalently the minimum number of powers of 2 needed to sum
to $m$ and let $v(0)=0$. 
(This is a simplified version of the example given in 
\cite[Example~9.17]{Da-La}.) 
Then $v$ is subadditive and unbounded on $\mathbb N_0$, but does not
tend to infinity as $m\to\infty$. 
The following lemma shows that the obvious extension of $v$ to $\mathbb R^+$
stays subadditive. Hence $\omega$ is an algebra weight with the
required properties.

\begin{lemma}
\label{le:subadd}
Let $v$ be a real-valued subadditive function on $\mathbb N_0$, and let also
$v$ denote the continuous extension of $v$ to $\mathbb R^+$ which is linear on each
of the intermediate intervals. Then $v$ be is subadditive on $\mathbb R^+$.
\end{lemma}

\noindent{\bf Proof}\quad
Let $x,y\in\mathbb R^+$ with $x=m+r$ and $y=n+s$, where $m,n\in\mathbb N_0$ and 
$0\le r,s<1$. Observe that $v(x)=(1-r)v(m)+rv(m+1)$ and similar for
$v(y)$. First, assume that $r+s\le1$. Then
\begin{align*}
v(x+y) &= (1-r-s)v(m+n)+(r+s)v(m+n+1)\\
&\le (1-r-s)(v(m)+v(n))+s(v(m)+v(n+1))+r(v(m+1)+v(n))\\
&=  (1-r)v(m)+(1-s)v(n)+sv(n+1))+rv(m+1)
=v(x)+v(y)
\end{align*}
as required. If $r+s>1$, we write $x+y=(m+n+1)+(r+s-1)$ and use the
same approach as above.
{\nopagebreak\hfill\raggedleft$\Box$\bigskip}


For homomorphisms $\Phi:L^1(\omega_1)\to L^1(\omega_2)$ several
conditions equivalent to the strong continuity of $(\nu^t)$ in $M(\omega_2)$ 
for $t\ge0$ is given in \cite[Corollary~3.13]{Gr} and 
\cite[Theorem~2.2]{Gh-Mc-Gr}. 
The proofs of most of these equivalencies carry over directly 
(with obvious modifications) to endomorphisms of $A(\omega)$, so by using
Theorem~\ref{th-sgcts}, we obtain the result below. We say that $f\in A(\omega)$ is a 
{\it standard element}, if the closed ideal it generates in $A(\omega)$ is
the {\it standard ideal} $A(\omega)_d=\{g\in A(\omega):\alpha(g)\ge d\}$, where 
$d=\alpha(f)=\inf\,\text{supp}(f)$. Also, we say that an endomorphism
$\Phi$ of $A(\omega)$ is a {\it standard endomorphism}, if whenever $f\in A(\omega)$
with $A(\omega)*f$ dense in $A(\omega)$, then $A(\omega)*\Phi(f)$ is dense in $A(\omega)$.

\begin{theorem}
\label{th-equiv}
Let $\Phi$ be a non-zero endomorphism of $A(\omega)$, let $(\nu^t)$ be the
semigroup from Theorem~\ref{th-hom} and let $(e_k)$ be a bounded
approximate identity for $A(\omega)$. Then the following are equivalent:
\begin{enumerate}[(a)]
\item 
$(\nu^t)$ is strongly continuous in $B(\omega)$ for $t\ge0$.
\item
There is a non-zero standard element $f\in A(\omega)$ for which 
$\nu^t*f\to f$ in $A(\omega)$ as $t\to0$.
\item
$\overline{\Phi(A(\omega))}$ contains a non-zero standard element.
\item
$\Phi$ is a standard endomorphism.
\item
$(\Phi(e_k))$ is a bounded approximate identity for $A(\omega)$.
\item
There is a non-zero standard element $f\in A(\omega)$ for which 
$\Phi(e_k)*f\to f$ in $A(\omega)$ as $k\to\infty$.
\end{enumerate}
In particular, these properties all hold if condition (\ref{eq:wein})
is satisfed.
\end{theorem}

\bigskip

Compared to the list in \cite[Theorem~2.2]{Gh-Mc-Gr} there are two
notable omissions in Theorem~\ref{th-equiv}, namely
\begin{enumerate}
\item[(g)]
For every $h\in A(\omega)$ there exist $f,g\in A(\omega)$ such that $h=\Phi(f)*g$.
\item[(h)]
The extension endomorphism $\widetilde{\Phi}$ of $B(\omega)$ is strongly continuous.
\end{enumerate}
Considering $A(\omega)$ as a Fréchet $A(\omega)$-module under the action 
$f\cdot g=\Phi(f)*g$ for $f,g\in A(\omega)$, property (g) can be restated as
existence of factorisation in this module, whereas (e) states that the
sequence $(e_k)$ is a bounded approximate identity for this module. 
As mentioned in Section~\ref{sec:various} we do not know whether
factorisation in $A(\omega)$ and its modules follows from the existence of
a bounded approximate identity 
(unless there exists $a>0$ such that 
$\sup_{n\in\mathbb N,\,0\le t\le a}\omega_n(t)<\infty$),  
that is, whether (e) implies (g).
However, the implications (g)$\Rightarrow$(h)$\Rightarrow$(a)-(f) follow
rather easily as in the proof of \cite[Theorem~2.2]{Gh-Mc-Gr}. 
Moreover, it may be possible to prove (h) from (a) (without assuming (g)).
For instance, it is easily seen that if $(\mu_{\beta})$ is a bounded
net in $B(\omega)$ with $\mu_{\beta}\to0$ strongly, then 
$\widetilde{\Phi}(\mu_{\beta})\to0$ strongly in $B(\omega)$.
We also mention in passing that if (a)-(f)
holds, then $\widetilde{\Phi}$ is $\sigma$-continuous 
(see Section~\ref{sec:various} for the definition): 
Let $(\mu_{\beta})$ be a net in $B(\omega)$ with 
$\mu_{\beta}\to0$ in the $\sigma$-topology, and let $n\in\mathbb N$. There exists
$m\in\mathbb N$ such that $\widetilde{\Phi}$ extends to a continuous homomorphism
$\widetilde{\Phi}_n:M(\omega_m)\to M(\omega_n)$.
Since $(\nu^t)$ is continuous in $B(\omega)$ and thus in $M(\omega_n)$ for $t\ge0$, it
follows from \cite[Theorem~2.2]{Gh-Mc-Gr} that $\widetilde{\Phi}_n$ is
strongly continuous. Since $\mu_{\beta}\to0$ strongly in $M(\omega_m)$ we thus have 
$\widetilde{\Phi}(\mu_{\beta})=\widetilde{\Phi}_n(\mu_{\beta})\to0$
strongly in $M(\omega_n)$. 
Hence $\widetilde{\Phi}(\mu_{\beta})\to0$ in the $\sigma$-topology on
$B(\omega)$, so $\widetilde{\Phi}$ is $\sigma$-continuous.

We end this section with a remark about automorphisms of $A(\omega)$.
Characterisations of automorphisms have been obtained for the
related algebras $L^1(\omega)$ (\cite{Gh:Hom} and \cite{Gh-Mc:Rad}), 
$L^1[0,1]$ (\cite{Ka-Sc}) and $L^1_{\text{loc}}(\mathbb R^+)$
(\cite{Gh-Mc}). For the Fr\'{e}chet algebra $L^1_{\text{loc}}(\mathbb R^+)$ 
(which is the projective limit of the algebras $L^1[0,n]\ (n\in\mathbb N)$) 
the method consists of reducing the problem to automorphisms of $L^1[0,1]$. 
The following simple example indicates, that questions about
automorphisms of $A(\omega)$ cannot easily be reduced to questions about
isomorphisms between the algebras $L^1(\omega_n)\ (n\in\mathbb N)$. Hence other
methods are needed in order to obtain a characterisation of the   
automorphisms of $A(\omega)$.
For $a\in\mathbb R^+$ we define an algebra weight $\omega_a$ on
$\mathbb R^+$ by 
$\omega_a(t)=e^{a\sqrt t}$ for $t\in\mathbb R^+$. Moreover, for a function $f$
on $\mathbb R^+$ we let $\Phi(f)(t)=2f(2t)$ for $t\in\mathbb R^+$. An
easy calculation
shows that $\|\Phi(f)\|_{\omega_a}=\|f\|_{\omega_{a/\sqrt 2}}$ for
$a\in\mathbb R^+$.
With $\omega=(\omega_n)$ this shows that $\Phi$ is an auotomorphism of
$A(\omega)$. However, for every $n\in\mathbb N$ there does not exist 
$m\in\mathbb N$ such
that $\Phi$ extends to an isomorphism between $L^1(\omega_m)$ and $L^1(\omega_n)$
(but $\Phi$ does extend to an isomorphism between 
$L^1(\omega_{n/\sqrt 2})$ and $L^1(\omega_n)$).

\section{Derivations}
\label{sec:der}

In Theorem~\ref{th:der} below we characterise the derivations on $A(\omega)$. 
Johnson (\cite{Jo:Cont} or \cite[Theorem~5.2.32]{Da:Book}) proved that
a semisimple, commutative Banach algebra does not have any non-zero 
derivations, and in particular this applies to each of the algebras $L^1(\omega_n)$.
For (radical) weights $\omega$, Jewell and Sinclair 
(see \cite[Remark~3(a)]{Je-Si} or \cite[Theorem~5.2.18\,(ii)]{Da:Book})  
proved that derivations on $L^1(\omega)$ are automatically continuous. Moreover, 
Ghahramani (\cite[Theorem~2.5]{Gh:Hom}) showed
that a linear operator $D$ on $L^1(\omega)$ is a derivation if and only if
there is a measure $\mu$ on $\mathbb R^+$ with
$$\sup_{t\in\mathbb R^+}\frac{t}{\omega(t)}\int\omega(t+s)\,d|\mu|(s)<\infty$$
such that $D(f)=(Xf)\ast\mu$ for $f\in L^1(\omega)$, where 
$(Xf)(t)=tf(t)$ for $t\in\mathbb R^+$ and $f\in L^1(\omega)$.  

Carpenter (\cite{Carp}; see also \cite[Theorem~8.2.5]{Go}) partly generalised
Johnson's result by showing that derivations on semisimple,
unital Fr\'{e}chet algebras automatically are continuous. We will see in
Corollary~\ref{co:autcts} that this implies that derivations on
$A(\omega)$ are continuous. (We remark that the proof of 
\cite[Theorem~5.2.18\,(ii)]{Da:Book} (using a gliding hump technique) 
cannot be used in this case, since $A(\omega)$ is not {\it locally bounded}, 
that is, has no bounded open sets 
(as a ball in $L^1(\omega_n)$ is unbounded in $L^1(\omega_m)$ for $m>n$).) 

Our main aim in this section is to prove the following result.

\begin{theorem}
\label{th:der}
\ 
\begin{enumerate}[(a)]
\item 
Suppose that 
\begin{equation}
\label{eq:weco}
\text{for every $n\in\mathbb N$ there exists $m\in\mathbb N$ such that}\quad
\sup_{t\in\mathbb R^+}\frac{t\omega_n(t)}{\omega_m(t)}<\infty.
\end{equation}
Then 
$$D_{\mu}(f)=(Xf)\ast\mu\qquad(f\in A(\omega))$$
defines a continuous derivation on $A(\omega)$ for every $\mu\in B(\omega)$ and
conversely every derivation on $A(\omega)$ has this form. 
\item
If condition (\ref{eq:weco}) is not satisfied, then there are no
non-zero derivations on $A(\omega)$.
\end{enumerate}
\end{theorem}

We remark that when condition (\ref{eq:weco}) is satisfied, then $X$ (and
thus $D_{\mu})$ extends to $B(\omega)$ by letting $d(X\nu)(t)=td\nu(t)$ for
$\nu\in B(\omega)$. Also, it is easily seem that condition
(\ref{eq:weco}) is equivalent with the formally weaker condition that
there exists some $p>0$ such that for every $n\in\mathbb N$ there exists
$m\in\mathbb N$ with 
$$\sup_{t\in\mathbb R^+}\frac{t^p\omega_n(t)}{\omega_m(t)}<\infty.$$

\medskip

Theorem~\ref{th:der} and its proof is closely related to Ghahramani's result.
Note however the difference
in that we only have a condition on the weights $(\omega_n)$ and that
if this condition is satisfied, then every $\mu\in B(\omega)$ gives rise to
a derivation.
We find it interesting that our results and methods for the ``little'' algebra
$A(\omega)$ are similar to those used for the ``big'' algebras 
$M_{\text{loc}}(\mathbb R^+)$ (\cite{Di}), $L^1_{\text{loc}}(\mathbb R^+)$ (\cite{Gh-Mc}) 
and $L^1(\omega)$ for radical weights $\omega$ (\cite{Gh:Hom}).

We need a few results in order to prove Theorem~\ref{th:der}. 
Recall that the identification in Theorem~\ref{th:mul} of $B(\omega)$ as
the multiplier algebra of $A(\omega)$ induces a strong operator topology on $B(\omega)$. 

\begin{lemma}
\label{le:derext}
Every derivation $D$ on $A(\omega)$ extends to a derivation $\overline{D}$
on $B(\omega)$ which is norm- as well as strongly continuous.
\end{lemma}

\noindent{\bf Proof}\quad
We proceed as in \cite[p.\,57]{Gh-Mc}: For $\mu\in B(\omega)$ we define a
linear map $S_{\mu}$ on $A(\omega)$ by 
$$S_{\mu}(f)=D(\mu*f)-\mu*D(f)\qquad(f\in A(\omega)).$$
It is easy to check that $S_{\mu}$ is a multiplier on $A(\omega)$. Hence it
follows from Theorem~\ref{th:mul} that there exists a unique measure 
$\overline{D}(\mu)\in B(\omega)$ such that $S_{\mu}(f)=\overline{D}(\mu)*f$
and thus 
$$D(\mu*f)=\overline{D}(\mu)*f+\mu*D(f)\qquad\text{for }f\in A(\omega).$$
It is further easily seen that the map $\mu\mapsto\overline{D}(\mu)$
defines a derivation on $B(\omega)$ which extends $D$.
Since $B(\omega)$ is semisimple by Corollary~\ref{co:ss}, it follows from
Carpenter's result mentioned above that $\overline{D}$ automatically
is continuous. 

For the strong continuity of $\overline{D}$, let
$(\mu_{\beta})$ be a net in $B(\omega)$ with $\mu_{\beta}\to0$ strongly. Then
$\overline{D}(\mu_{\beta})*f=D(\mu_{\beta}*f)-\mu_{\beta}*D(f)\to0$
in $A(\omega)$ for every $f\in A(\omega)$, so 
$\overline{D}(\mu_{\beta})\to0$ strongly in $B(\omega)$ as required.
{\nopagebreak\hfill\raggedleft$\Box$\bigskip}

\begin{corollary}
\label{co:autcts}
Derivations on $A(\omega)$ are continuous.
\end{corollary}

The next result is similar to \cite[Lemma~1]{Di}.
Since the notation
is somewhat different we include a proof for the sake of completeness.
For a measure $\mu$ on $\mathbb R^+$ we let $\alpha(\mu)=\inf\,\text{supp}(\mu)$.

\begin{lemma}
\label{le:dersupp}
For a derivation $D$ on $B(\omega)$ we have
$$\alpha(D(\mu))\ge\alpha(\mu)\qquad\text{for }\mu\in B(\omega).$$
In particular $\alpha(D(\delta_t))\ge t$ for $t\in\mathbb R^+$.
\end{lemma}

\noindent{\bf Proof}\quad
Let $t\in\mathbb R^+$ and $k\in\mathbb N$. Since $\delta_t=(\delta_{t/k})^{*k}$ we
have $D(\delta_t)=k(\delta_{t/k})^{*(k-1)}*D(\delta_{t/k})$, so it
follows from Titchmarsh's convolution theorem that 
$$\alpha(D(\delta_t))\ge\alpha((\delta_{t/k})^{*(k-1)})=\tfrac{(k-1)t}{k}.$$
Hence $\alpha(D(\delta_t))\ge t$. Now let $\mu\in B(\omega)$ and let
$t=\alpha(\mu)$. We may assume that $t>0$. Define the translate
$\mu_t\in B(\omega)$ by $\mu_t(E)=\mu(E+t)$ for measurable sets
$E\subseteq\mathbb R^+$. Then $\mu=\delta_t*\mu_t$, so 
$D(\mu)=D(\delta_t)*\mu_t+\delta_t*D(\mu_t)$, and thus
$$\alpha(D(\mu))\ge\min\{\alpha(D(\delta_t)*\mu_t),\alpha(\delta_t*D(\mu_t))\}
\ge\min\{\alpha(D(\delta_t)),\alpha(\delta_t)\}=t$$
as required.
{\nopagebreak\hfill\raggedleft$\Box$\bigskip}

The proof of the following result is identical to that of 
\cite[Lemma~2.3]{Gh:Hom} (see also \cite[Lemma~3]{Di}) and is
therefore omitted.

\begin{lemma}
\label{le:dermu}
Let $D$ be a derivation on $B(\omega)$. Then there exists $\mu\in B(\omega)$ such
that
$$D(\delta_t)=t\cdot\delta_t*\mu\qquad\text{for }t\in\mathbb R^+.$$
\end{lemma}

\bigskip

Proposition~\ref{prop:dualSOsi} plays a crucial role in the proof of 
the next result.

\begin{lemma}
\label{le:dense}
The linear span of $\{\delta_t:t\in\mathbb R^+\}$ is strongly dense in $B(\omega)$.
\end{lemma}

\noindent{\bf Proof}\quad
Let $V$ be the linear span of $\{\delta_t:t\in\mathbb R^+\}$ and let $\varphi$
be a strongly continuous linear functional on $B(\omega)$ with $\varphi(\delta_t)=0$
for $t\in\mathbb R^+$. From the proof of Proposition~\ref{prop:dualSOsi} it
follows that $\varphi$ extends to a strongly continuous linear functional 
$\varphi_N:M({\omega_N})\to\mathbb C$ for some $N\in\mathbb N$. However, $V$ is 
strongly dense in $M({\omega_N})$ by \cite[Lemma~1.3]{Gh:Hom} 
(see also \cite[Lemma~4]{Di}), so we deduce that $\varphi_N$ and thus
$\varphi$ is zero, which finishes the proof.
{\nopagebreak\hfill\raggedleft$\Box$\bigskip}

We are now ready to prove the characterisation of derivations on $A(\omega)$.
\medskip

\noindent{\bf Proof of Theorem~\ref{th:der}}\quad
\begin{enumerate}[(a)]
\item 
If condition (\ref{eq:weco}) is satisfied, then $X$ defines a continuous
linear operator $X:L^1(\omega_m)\to L^1(\omega_n)$ for every $n\in\mathbb N$ with corresponding
$m\in\mathbb N$, and thus defines a continuous linear operator on $A(\omega)$. A
routine calculation shows that $X$ is a derivation. Moreover, for
$\mu\in B(\omega)$ the map $T_{\mu}(f)=\mu\ast f\ (f\in A(\omega))$ defines a continuous 
multiplier on $A(\omega)$ by Theorem~\ref{th:mul}. Hence $D_{\mu}=T_{\mu}X$
is a continuous derivation on $A(\omega)$.

Conversely, let $D$ be a derivation on $A(\omega)$ and let 
$\overline{D}$ be the continuous extension to a derivation on $B(\omega)$ given by
Lemma~\ref{le:derext}. Also, let $\mu\in B(\omega)$ be the measure from 
Lemma~\ref{le:dermu} satisfying 
$\overline{D}(\delta_t)=t\cdot\delta_t*\mu$ for $t\in\mathbb R^+$.
Define a continuous derivation $\Delta$ on $B(\omega)$ by
$$\Delta(\nu)=\overline{D}(\nu)-(X\nu)*\mu\qquad\text{for }\nu\in B(\omega).$$
Then $\Delta(\delta_t)=0$ for $t\in\mathbb R^+$. Furthermore, it follows from
the proof of Lemma~\ref{le:derext} that $\Delta$ is strongly continuous. Since 
the linear span of $\{\delta_t:t\in\mathbb R^+\}$ is strongly dense in $B(\omega)$ by
Lemma~\ref{le:dense}, we thus deduce that $\Delta=0$, so
$\overline{D}(\nu)=(X\nu)*\mu$ for $\nu\in B(\omega)$ and thus $D=D_{\mu}$
as required.
\item
Assume that $D$ is a non-zero derivation on $A(\omega)$. As
in the proof of (a), $D$ extends to a continuous derivation $\overline{D}$
on $B(\omega)$ and there exists $\mu\in B(\omega)$ such that 
$\overline{D}(\delta_t)=t\cdot\delta_t*\mu$ for $t\in\mathbb R^+$. 
Given $n\in\mathbb N$ it follows from Lemma~\ref{le:top} that there exists
$m\in\mathbb N$ such that $D$ extends continuously to a map $D:L^1(\omega_m)\to L^1(\omega_n)$. 
Hence Grabiner's extension (\cite[Theorem~3.4]{Gr}) 
$\overline{D}:M(\omega_m)\to M(\omega_n)$ is also an extension of the derivation 
$\overline{D}$ on $B(\omega)$, so there exists a constant $c_n$ such that
$\|\overline{D}(\delta_t)\|_n\le c_n\|\delta_t\|_m$ for all $t\in\mathbb R^+$. 
However, 
$\|\overline{D}(\delta_t)\|_n=\int_0^{\infty} t\omega_n(s+t)\,d|\mu|(s)
\ge t\omega_n(t)|\mu|(\mathbb R^+)$, whereas $\|\delta_t\|_m=\omega_m(t)$, so 
we conclude that (\ref{eq:weco}) is satisfied.
\end{enumerate}
{\nopagebreak\hfill\raggedleft$\Box$\bigskip}

\section{Summary of open problems}

We finish by gathering some open problems that have been mentioned
in this paper. 
\begin{enumerate}[(I)]
\item
Do we have factorisation in $A(\omega)$ (and in Fr\'echet $A(\omega)$-modules) if 
$\lim_{n\to\infty}\omega_n(t)=\infty$ for every $t\in\mathbb R^+$? 
(See the discussion preceeding Theorem~\ref{th:mul}.)
\item
Are the $\sigma$- and the SO-topologies identical?
(See the end of Section~\ref{sec:various}.) 
\item
Are the conditions 
\begin{enumerate}
\item[(g)]
For every $h\in A(\omega)$ there exist $f,g\in A(\omega)$ such that $h=\Phi(f)*g$.
\item[(h)]
The extension endomorphism $\widetilde{\Phi}$ of $B(\omega)$ is strongly continuous.
\end{enumerate}
equivalent to the conditions in Theorem~\ref{th-equiv}?
(See the discussion following Theorem~\ref{th-equiv}.)
\end{enumerate}

\bigskip

\noindent
Thomas Vils Pedersen\\
Department of Basic Sciences and Environment\\
Faculty of Life Sciences\\
University of Copenhagen\\
Thorvaldsensvej 40\\
DK-1871 Frederiksberg C\\
Denmark\\
vils@life.ku.dk

\end{document}